\documentclass{birkjour}
\usepackage{amssymb}
\usepackage{amsthm,cancel}

\usepackage[colorlinks]{hyperref}
\usepackage{cleveref}
\usepackage{color,graphicx,shortvrb}

\usepackage[active]{srcltx} 

\usepackage{enumerate}
\usepackage{amssymb}
\usepackage{tikz-cd}

\usepackage{mathtools}

 \newtheorem{theorem}{Theorem}[section]
 \newtheorem{corollary}{Corollary}[section]
 \newtheorem{lemma}{Lemma}[section]
 \newtheorem{proposition}{Proposition}[section]
 \theoremstyle{definition}
 
 \theoremstyle{remark}
 \newtheorem{remark}{Remark}[section]
 \newtheorem{example}{Example}
 
 \numberwithin{equation}{section}
 \usepackage{amsmath, amssymb, amsthm}
 \usepackage{mathrsfs}

\begin{document}

%
%
%
%
%
%
%
%
%

\title{New insights into linear maps {which are} anti-derivable at zero}

\author[J. Li]{Jiankui Li}
\address{%
School of Mathematics, East China University of Science and Technology, Shanghai, 200237 China.}
\email{jkli@ecust.edu.cn}

\author[A.M. Peralta]{Antonio M. Peralta}
\address{Instituto de Matem{\'a}ticas de la Universidad de Granada (IMAG), Departamento de An{\'a}lisis Matem{\'a}tico, Facultad de
	Ciencias, Universidad de Granada, 18071 Granada, Spain.}
\email{aperalta@ugr.es}

\author[S. Su]{Shanshan Su}
\address{%
(Current address) School of Mathematics, East China University of Science and Technology, Shanghai, 200237 China. \\
Departamento de An{\'a}lisis Matem{\'a}tico, Facultad de
	Ciencias, Universidad de Granada, 18071 Granada, Spain.}
\email{lat875rina@gmail.com}




\date{today}

\begin{abstract} Let $A$ be a Banach algebra admitting a bounded approximate unit and satisfying property $\mathbb{B}$. Suppose $T: A \rightarrow X$ is a continuous linear map, where $X$ is an essential Banach $A$-bimodule. We prove that the following statements are equivalent:
\begin{enumerate}[$(i)$]
\item $T$ is anti-derivable at zero (i.e., $a b =0$ in $A$ $\Rightarrow T(b)\cdot a + b\cdot T(a) =0$);
\item There exist an element $\xi \in  X^{**}$ and a linear map (actually a bounded Jordan derivation) $d: A\to X$ satisfying $\xi \cdot a = a \cdot \xi \in X$, 
$T(a) = d(a) +\xi \cdot a$, and $d(b)\cdot a + b\cdot d(a)=  - 2 \xi \cdot (b a),$ for all $a,b\in A$ with $a b =0$. 
\end{enumerate} 
Assuming that $A$ is a C$^*$-algebra we show that a bounded linear mapping $T: A\to X$ is anti-derivable at zero if, and only if, there exist an element $\eta \in  X^{**}$ and an anti-derivation $d: A \rightarrow X$ satisfying $\eta \cdot a = a \cdot \eta \in X$, $\eta \cdot [a,b] = 0$ {\rm(}i.e., $L_{\eta}: A \to A$, $L_{\eta} (a) = \eta \cdot a$ vanishes on commutators{\rm)}, and $T(a) = d(a) +\eta \cdot a$, for all $a,b \in A$. The results are also applied for some special operator algebras.
\end{abstract}

\maketitle
\section{Introduction}\label{intro}

It is well known that linear homomorphisms and anti-homomorphisms  (i.e., linear maps preserving {and} reversing products, {respectively}) between associative algebras both preserve squares of elements, and hence Jordan products of the form $a\circ b = \frac12 (a b + ba )$. Actually, direct sums of homomorphisms and anti-homomorphisms are enough to describe all Jordan $^*$-homomorphisms from a C$^*$-algebra $A$ into the von Neumann algebra, $B(H),$ of all bounded operators on a complex Hilbert space $H$ (cf. \cite[Theorem 10]{Kad51} and \cite[Theorem 3.3]{Sotr65}), and similarly for Jordan homomorphisms from a ring onto a 2-torsion free semiprime ring (cf. \cite[Theorem 2.3]{Bre1989}). It is quite natural to ask whether anti-derivations could have such a central role in the case of Jordan derivations. Let us recall that a linear mapping $d$ from an associative algebra $A$ into an $A$-bimodule $X$ is called a {\emph{derivation}} (respectively, an {\emph{anti-derivation} or a \emph{Jordan derivation}}) if $d (a b ) = d(a)\cdot b + a\cdot d(b)$ (respectively, $d (a b ) = d(b)\cdot a + b\cdot d(a),$ or $d(a\circ b)= d(a) \circ b + a \circ d(b)$) for all $a,b\in A$. In case that $A$ is commutative, derivations and anti-derivations {define} the same maps. \smallskip

A milestone result due to B.E. Johnson proves that every bounded Jordan derivation from a C$^*$-algebra $A$ to a Banach $A$-bimodule is a derivation (cf. \cite{Johnson2001localderivationsonC}). Recall that derivations, Jordan derivations, and hence anti-derivations from a C$^*$-algebra $A$ into a Banach $A$-bimodule are all continuous (\cite{Ringrose72} and \cite[Corollary 17]{PeRusso2014automaticcontinuity}). So, every Jordan derivation from a C$^*$-algebra $A$ into a Banach $A$-bimodule is a derivation. It was recently established that anti-derivations whose domain is a C$^*$-algebra are quite strange maps, namely, if $\delta: A\to X$ is an anti-derivation from a C$^*$-algebra to a Banach $A$-bimodule, then there exists a finite central projection $p_1$ in $A^{**}$ and an element $x_0\in X^{**}$ such that $$\delta^{**} (x) = \delta_{x_0} (\tau(p_1 x)) = [\tau(p_1 x) , x_0], \hbox{ for all } x\in A^{**},$$ where $\tau: p_1 A^{**}\to Z(p_1 A^{**})$ is the (faithful) center-valued trace on $p_1 A^{**}$ (see \cite[Theorem 4]{abulhamil2020linear}). Nevertheless, the just quoted reference shows the existence of quite natural examples of non-zero anti-derivations. Newly results, due to W. Huang, S. Li and the first author of this note, establish that every anti-derivation from a C$^*$-algebra $A$ acting on a Hilbert space $H$ into $B(H)$ is inner, furthermore, if the commutant of all commutators in $A$ lies in the commutant of $A$ in $B(H)$, every anti-derivation from $A$ into $B(H)$ is zero \cite{HuangLis2025antiderivations}. The just quoted reference also contains results showing that for each semiprime Banach algebra $A$ with semisimple center, every bounded {anti-derivation} on $A$ is zero, and in case that $A$ is a semisimple Banach algebra with nonzero socle, there are no nonzero {anti-derivations}
from the socle of $A$ into $A$. The same authors also establish conditions to guarantee the existence of non-trivial anti-derivations on nest algebras.\smallskip  

The current research has devoted some efforts to determining anti-derivations in different classes of algebras, as well as to describing those maps satisfying the weaker property of being anti-derivable at zero. A linear map $T$ from an algebra $A$ into an $A$-bimodule $X$ is called \emph{anti-derivable at zero} 
$$d(b) \cdot a + b \cdot d(a) = 0 \text{ whenever }a b  =0 \text{ in } A.$$ B. Fadaee and H. Ghahramani proved in \cite[Theorem 3.3]{FadaeeHoger2020linear} that a bounded linear map $T$ from a C$^*$-algebra $A$ into $A^{**}$ is anti-derivable at zero if, and only if, there exist a continuous derivation $d : A \to A^{**}$ and a central element $\eta\in Z(A^{**})$ satisfying $T(a) = d (a) + \eta \cdot a$ for all $a \in  A$. If $X$ is a complex Banach space with dim$(X)\geq 2$, and $A\subseteq B(X)$ is a standard operator algebra, the unique linear mapping on $A$ which is anti-derivable at zero is the zero mapping, a result due to K. Fallahi and H. Ghahramani \cite{FallahiHoger2022anti}.\smallskip 

D.A. Abulhamil, F.B. Jamjoom, and the second author of this note established in \cite[Theorem 6]{abulhamil2020linear} the first result on bounded linear maps from a C$^*$-algebra $A$ into an essential Banach $A$-bimodule $X$ which are anti-derivable at zero. Concretely, a bounded linear mapping $T: A\to X$ is anti-derivable at zero if, and only if, there exist $\eta \in  X^{**}$ and an anti-derivation $d: A \rightarrow X^{**}$ satisfying $\eta \cdot a = a \cdot \eta \in X$, $\eta \cdot [a,b] = 0$, and $T(a) = d(a) +\eta \cdot a$, for all $a,b \in A$.\smallskip 

More recently, L. Liu and S. Hou have completely determined all linear maps on certain generalized matrix algebras which are anti-derivable at zero, showing that these maps are expressed as the sum of an anti-derivation $d_1$, a derivation $d_2$, and the multiplication operator by a central element $\eta$ satisfying $d_2 ([a, b]) = - 2 \eta [a, b]$ for all $a,b$ with $a b =0$ (cf. \cite{Liu2025characterization}).\smallskip

In this note we shall resume the study of bounded linear maps anti-derivable at zero whose domain is a C$^*$-algebra. Firstly, by improving the conclusion in \cite[Theorem 6]{abulhamil2020linear} above, by showing that we can take an anti-derivation $d: A\to X$. We must additionally point out that we have detected a gap in the proof originally given in \cite{abulhamil2020linear}. The problem arises in assuming that for each Banach $A$-bimodule $X$ over a C$^*$-algebra $A$, the opposite module, $X^{op},$ is a Banach $A$-bimodule too, an statement that may {fail} in some cases (cf. the comments before \Cref{thm C-star algebra being anti-derivable at zero}). As we shall see later, here we also provide a complete new proof of this result which avoids the commented difficulties; so, the problem only affects to the arguments but not to the conclusion in \cite[Theorem 6]{abulhamil2020linear}, which remains completely valid. Section~\ref{sec C-star-algebras} is entirely devoted to this task. We also include a couple of corollaries for C$^*$-algebras whose second dual does not admit type $I_1$ summands and for von Neumann algebras, where the continuity of the linear mapping can be relaxed.\smallskip 

In addition to  the already commented fact that every Jordan derivation from a C$^*$-algebra $A$ to a Banach $A$-bimodule is a derivation, each C$^*$-algebra possesses a bounded approximate unit (\cite[Theorem 3.2.21]{Dales2001banachalgs}) and satisfies property $\mathbb{B}$ (\cite{ABEV2009mapspreservingzpd} and \cite[Theorem 5.19]{Brevsar2021zpdalgs}). We begin our study with the case of  bounded linear maps whose domain is a Banach algebra satisfying the above properties. We establish that if $A$ is a Banach algebra satisfying property $\mathbb{B}$ and admitting a bounded approximate unit, a bounded linear operator $T$ from $A$ into a Banach $A$-bimodule $X$ is anti-derivable at zero if, and only if, there {exist} an element $\xi \in  X^{**},$ and a linear map (actually a bounded Jordan derivation) $d: A\to X$ satisfying $\xi \cdot a = a \cdot \xi \in X$, 
$T(a) = d(a) +\xi \cdot a$, and $d(b)\cdot a + b\cdot d(a)=  - 2 \xi \cdot (b a),$ for all $a,b\in A$ with $a b =0$. If we also assume that  every Jordan derivation from $A$ to $X$ is a derivation, we can assume that $d$ above is a derivation  satisfying $d([a,b]) = -2\xi\cdot [a,b]$, for all $a,b\in A$ (see \Cref{thm anti-derivable at zero}). It should be noted that in case that $A$ and $X$ are {both unital, and} $A$ satisfies property $\mathbb{B}$, the problem of determining those linear maps from $A$ into $X$ which are anti-derivable at zero admits a simpler solution without assuming continuity of these maps (cf. \Cref{r the unital case in Thm 3.1 is trivial}). Applications of our main result are found in the cases of C$^*$-algebras, algebras of the form $M_n (R)$, where $R$ is a {2-torsion free} unital (associative) ring and $n\geq 2$, and nest algebras. \smallskip


\subsection{Technicalities}

Unless otherwise stated, all the algebras in this paper will be over the complex field. Zero product determined algebras constitute a rich class of algebras satisfying interesting properties which can be applied in different problems. According to \cite{Brevsar2021zpdalgs}, we say that an algebra $A$ is \emph{zero product determined} if for every linear space $X$, and every bilinear map $\varphi: A \times A \to X$ satisfying 
$$\varphi(a,b) = 0, \hbox{ for all } a,b \in A \text{ with } a b  =0 ,$$
there exists a linear map $G: \hbox{span}(\{ a b : a, b\in A\}) \to X$ such that 
$$\varphi(a, b ) = G(a b),$$ for all $a,b \in A$. If $A$ is a Banach algebra and $X$ is a Banach space, then we also require that $\varphi$ and $G$ are continuous. A (Banach) algebra $A$ is said to satisfy \emph{property $\mathbb{B}$} (see \cite{ABEV2009mapspreservingzpd}) if for every (Banach) $A$-bimodule $X$, every (continuous) bilinear map $\varphi: A \times A \to X$ satisfying 
$$\varphi(a,b) = 0, \text{ for all } a,b \in A \text{ with } a b =0 ,$$ also satisfies the identity 
$$\varphi(a b , c) = \varphi(a, b c),  \text{ for all } a,b,c \in A.$$

It is clear that $A$ being zero product determined implies that $A$ has property $\mathbb{B}$. However, the reciprocal implication is not, in general, true (see \cite[Example 5.3]{Brevsar2021zpdalgs}). Lemma 2.3 in \cite{ABEV2009mapspreservingzpd} assures that if $A$ is a Banach algebra admitting a bounded left approximate identity, then $A$ is zero product determined if, and only if, $A$ has property $\mathbb{B}$.
\smallskip

Examples of Banach algebras satisfying property $\mathbb{B}$ include group algebras $L^1(G)$ for any locally compact group $G$, C$^*$-algebras, the Banach algebras of all approximable operators and all nuclear operators on any Banach space $X$, the {$p$-Schatten von} Neumann classes $S_p(H)$, for any Hilbert space $H$ and any {$1 \leq 
p \leq \infty$}, the Banach algebras $\ell^p(I)$ with $1 \leq  p < \infty$ and $c_0(I)$ for any nonempty set $I$, among others {(see, for instance, \cite{ABEV2009mapspreservingzpd}, \cite[\S 5]{Brevsar2021zpdalgs} and \cite[Theorem 2.11 and Examples in 1.3]{Bresar2009zpdmatrixalgebras})}.\smallskip

Along this paper, the center of an algebra $A$, that is, the set of all elements in $A$ that commute with every element in $A$, will be denoted by $\mathcal{Z}(A)$.\smallskip
    
It is well known that the bidual of $A$ ($A^{**}$ from now on) admits at least two different Arens products making it a Banach algebra (see \cite[Definition 2.6.16]{Dales2001banachalgs} or \cite[\S 1.4]{PalmerBook1994}). We shall focus on the first Arens product. Furthermore, the bidual, $X^{**}$, of any Banach {$A$-bimodule $X$}, can be naturally endowed with a Banach $A^{**}$-bimodule structure through the first Arens product (see the construction in \cite[pages 248 and 249]{Dales2001banachalgs}). The separate weak$^*$-continuity properties of the Arens product can be consulted in \cite[Proposition A.3.52]{Dales2001banachalgs} and \cite[page 48]{PalmerBook1994}. We recall a property employed in later arguments: if $(a_{\lambda})_{\lambda}$ and $(x_{\mu})_{\mu}$ are nets in $A$ and $X$, respectively, such that $(a_{\lambda})_{\lambda}$ converges to $a \in A^{**}$ in the weak$^*$-topology of $A^{**}$ and $(x_{\mu})_{\mu}$ converges to $x \in X^{**}$ in the weak$^*$-topology of $X^{**}$, then 
\begin{equation}\label{eq weak*-continuity properties} a \cdot x = w^{*}\text{-} \lim_{\lambda} \lim_{\mu} a_{\lambda} \cdot x_{\mu} \text{ and } x \cdot a = w^{*}\text{-}\lim_{\mu} \lim_{\lambda} x_{\mu} \cdot a_{\lambda},
\end{equation}
in $X^{**}$ (see \cite[(2.6.26)]{Dales2001banachalgs}).\smallskip 

Under the above conditions, given a bounded left approximate unit $(e_j)_{j}$ in $A$, if $e$ denotes any weak$^*$ cluster point of the net $(e_j)_{j}$ in $A^{**}$. Then $e\cdot a = a$
for every element $a \in A$. Having in mind that the first Arens product is weak$^*$ continuous when we fix an arbitrary element in $A^{**}$ in the second variable, we obtain that $e a = a$ for all $a\in A^{**}$. The converse also holds (cf. \cite[Proposition 5.1.8]{PalmerBook1994}). Note that if $A$ is Arens regular, we can actually deduce that $a e = a$, for all $a\in A^{**}$.\smallskip

\begin{remark}\label{r essential modules and bau} A Banach $A$-bimodule $X$ over a Banach algebra $A$ is called \emph{essential} if the linear span of the set $\{a \cdot x \cdot b: a,b \in A, x \in X\}$ is dense in $X$. Suppose now that $X$ is essential, and $A$ admits a bounded approximate unit $(e_j)_j$. As before, let $e$ be a weak$^*$ cluster point of $(e_j)_j$ in $A^{**}$. It is not hard to see that, by considering $X$ as a subspace of $X^{**}$ which is a Banach $A^{**}$-bimodule under the first Arens product, we have $\xi \cdot e = e \cdot \xi = \xi,$ for every $\xi \in X.$ Consequently, $\xi \cdot e = \xi,$ for all $\xi\in X^{**}$.\smallskip	
	
Furthermore, if we combine \eqref{eq weak*-continuity properties} with the properties of the bounded approximate unit, it can be easily seen that $e^2 = e$ is a projection in $A^{**}$. 	
\end{remark}

\section{Linear maps which are anti-derivable at zero whose domain is a zero product determined Banach algebra}

The main goal of this section is to characterize continuous linear maps anti-derivable at zero when the domain is a Banach algebra with some reasonable assumptions. We begin with some examples to observe that non-trivial anti-derivations and linear maps being anti-derivable at zero do exist.

\begin{example}\label{ex nontrivial anti-derivations}
 Let $A= \left\{\begin{pmatrix}
        a & b \\ 
        0 & c
    \end{pmatrix}: a,b,c \in \mathbb{C} \right\}$ denote the algebra of all $2$ by $2$ upper triangular {matrices} over $\mathbb{C}$. Fix any $\lambda \in \mathbb{C}$, and define a linear map $d_{\lambda}: A \to A$ by $d\begin{pmatrix}
        a & b \\ 
        0 & c
    \end{pmatrix} = \begin{pmatrix}
        0 & \lambda(a-c) \\ 
        0 & 0
    \end{pmatrix}$. It is not hard {to check} that $d_{\lambda}$ is an anti-derivation on $A$.
\end{example}

\begin{example}\label{ex nontrivial anti-derivable at zero}
 Consider a linear map $T:A \to A$ defined by $T\begin{pmatrix}
        a & b \\ 
        0 & c
    \end{pmatrix} = \begin{pmatrix}
        a & -b \\ 
        0 & c
    \end{pmatrix}$, where $A$ is the algebra in the previous example. Let $\mathbf{1}$ denote the unit of $A$. It is easy to see that $T$ is anti-derivable at zero. Moreover, the {map} $d: A \to A $ defined by $d\begin{pmatrix}
        a & b \\ 
        0 & c
    \end{pmatrix} = \begin{pmatrix}
        0 & -2 b \\ 
        0 & 0
    \end{pmatrix}$ is a derivation, and thus a Jordan derivation on $A$, which is not an anti-derivation, and the identity $T(x) = d(x) + \mathbf{1} x,$ holds for all $x\in A$. Furthermore, if $x y =0$ in $A$, we have $d(y) x + y d(x) = -2 y x= 2 [x,y].$\smallskip 
    
We claim that it is impossible to find a non-zero anti-derivation $\tilde{d}$ on $A$ and $\xi \in \mathcal{Z}(A)$ such that $T(x) = \tilde{d}(x) + \xi x,$ for all $x\in A$. Namely, since $\mathcal{Z}(A) = \mathbb{C} \mathbf{1},$ if there is a non-zero anti-derivation $\tilde{d}: A \to A$ such that $T(x) = \tilde{d}(x) + \xi x = \tilde{d}(x) + \alpha x$ for some $\alpha \in \mathbb{C}$, then $$\tilde{d}(x) = T(x) - \alpha x, \hbox{ and } \tilde{d}(xy) = \tilde{d}(y) x + y \tilde{d}(x),$$ for all $x,y \in A$. However,  by considering, for instance, $x = \begin{pmatrix}
        1 & 0 \\ 
        0 & 0
    \end{pmatrix}$ and $y =\begin{pmatrix}
        1 & 1 \\ 
        0 & 1
    \end{pmatrix}$, we have $$\tilde{d} (xy) = \begin{pmatrix}
    1-\alpha & -1-\alpha \\ 
    0 & 0
    \end{pmatrix}, \text{ and } \tilde{d}(y)x + y \tilde{d}(x) = \begin{pmatrix}
        2-2\alpha & 0 \\ 
        0 & 0
    \end{pmatrix}$$
    which is unsolvable for $\alpha\in \mathbb{C}$.
\end{example}

According to the most used notation (see, for example, \cite{ABEV2009mapspreservingzpd,HPS2025}), a linear map $T$ from {an associative} Banach algebra $A$ into a Banach $A$-bimodule $X$ is called a \emph{generalized derivation} if there exists $\xi\in X^{**}$ satisfying $$T(x y ) = T(x) y + x T(y)  - x\cdot \xi \cdot y, \hbox{ for all } x,y\in A.$$ A linear map $\delta: A\to X$ is called a \emph{Jordan derivation} if $\delta (a^2) = \delta (a) \cdot a + a\cdot \delta (a),$ for all $a\in A$. Clearly, every derivation is a Jordan derivation, and the problem of determining under which conditions every Jordan derivation is actually a derivation has attracted a lot of attention (cf. \cite{Johnson1996,Lu2009JordanCSLalgebras}).\smallskip

The main result in this section characterises continuous linear maps which are anti-derivable at zero whose domain is a Banach algebra $A$ with property $\mathbb{B}$ (cf. \cite[Theorem 3.1]{Liu2025characterization} for a result on certain generalized matrix algebras). Note that in view of \Cref{ex nontrivial anti-derivable at zero}, a bounded linear map anti-derivable at zero need not be necessarily written as the sum of an anti-derivation and the left multiplication operator by an element in the centre.

\begin{theorem}\label{thm anti-derivable at zero} Let $A$ be a Banach algebra satisfying property $\mathbb{B}$ and admitting a bounded approximate unit, and let $X$ be an essential Banach $A$-bimodule. Then the following statements are equivalent for every continuous linear map $T: A \rightarrow X$:
	\begin{enumerate}[$(i)$]
		\item $T$ is anti-derivable at zero.
		\item There {exist} an element $\xi \in  X^{**},$ and a linear map (actually a bounded Jordan derivation) $d: A\to X$ satisfying $\xi \cdot a = a \cdot \xi \in X$, 
		$T(a) = d(a) +\xi \cdot a$, and $d(b)\cdot a + b\cdot d(a)=  - 2 \xi \cdot (b a),$ for all $a,b\in A$ with $a b =0$. 
		\end{enumerate}
If we additionally assume that every Jordan derivation from $A$ to $X$ is a derivation, statements $(i)$ and $(ii)$ above are equivalent to:	\begin{enumerate}[$(iii)$]
	\item There exists an element $\xi \in  X^{**},$ and a derivation $d: A\to X$ satisfying $\xi \cdot a = a \cdot \xi \in X$, $T(a) = d(a) +\xi \cdot a$, and $d([a,b]) = -2\xi\cdot [a,b]$, for all $a,b\in A$. In particular, $T$ is a generalized derivation.
\end{enumerate}
\end{theorem}

\begin{proof}$(i) \Rightarrow (ii)$ Assume that $T$ is anti-derivable at zero. Let $\varphi : A \times A \rightarrow X$ be the bilinear map defined by $\varphi(a,b) = T(b) \cdot a + b\cdot T(a)$. Then, by hypotheses, $\varphi$ is continuous and $a b = 0$ in $A$ implies $\varphi(a,b) = 0$. Since $A$ satisfies property $\mathbb{B}$, the identity
\begin{equation}\label{eq key identity 2209} T(c) \cdot (a b) + c\cdot T(a b) = \varphi(a b,c) = \varphi(a,b c) = T(b c) \cdot a + (b c)\cdot T(a),
\end{equation} holds for all $a,b,c \in A$. Let $(e_j)_j$ denote a bounded approximate unit in $A$, and let $e$ be a weak$^*$ cluster point of $(e_j)_j$ in $A^{**}$. By replacing {$c$} with $e_j$, and taking weak$^*$ and norm limits of a suitable subnet we arrive, via \Cref{r essential modules and bau} and \eqref{eq weak*-continuity properties}, to \begin{equation}\label{eq 2 on 2209} T^{**}(e) \cdot (a b) + T(a b) = T^{**}(e) \cdot (a b ) + e\cdot T(a b) = T(b) \cdot a + b\cdot T(a),
\end{equation} for all $a,b\in A$.  By replacing $a$ with $e_{j},$ and taking weak$^*$ and norm limits of an appropriate subnet in the previous identity, we get 
$$T^{**}(e) \cdot b + T(b) =  T(b) \cdot e + b\cdot T^{**}(e) =  T(b) + b\cdot T^{**}(e),$$ for all $b\in A$. That is, $\xi = T^{**}(e)\in X^{**}$ commutes with all elements in $A$. By combining this information with the identity in \eqref{eq 2 on 2209} it easily follows that \begin{equation}\label{eq T is a kind of generalized antider} \begin{aligned}
T(a b) = T(b) \cdot a + b\cdot T(a)- a\cdot  \xi \cdot b &= T(b) \cdot a + b\cdot T(a)- \xi \cdot (a b) \\ 
&= T(b) \cdot a + b\cdot T(a)-(a b) \cdot \xi, 
	\end{aligned}
\end{equation} for all $a,b\in A$. 
\smallskip 

Define $d: A \rightarrow X$ by $d(a) = T(a) - \xi \cdot a$ {for all $a$ in $A$}. Having in mind \eqref{eq T is a kind of generalized antider}, it can be easily deduced that $$\begin{aligned}
d(a^{2}) &= T(a^{2}) - \xi \cdot a^{2} = T(a) \cdot a + a\cdot T(a) - 2\xi \cdot a^{2} =  d(a) \cdot a + a \cdot d(a),
\end{aligned}$$ which confirms that $d$ is a continuous Jordan derivation, and clearly $T(a) = d(a) + \xi \cdot a,$ for all $a\in A$.\smallskip 

Observe that if $a b =0 $ in {$A$,} we must have
$$ 0 = T(b)\cdot a + b \cdot T(a)  = d(b) \cdot a + \xi \cdot (b a) + b \cdot d(a) + \xi \cdot (b a ), $$ which assures that $ d(b) \cdot a + b \cdot d(a) = - 2 \xi \cdot (b a)$, and gives the desired statement in $(ii)$.\smallskip

$(ii) \Rightarrow (i)$ Suppose that we can write $T(a) = d(a) + \xi\cdot a$ for all $a\in A$, where {$\xi\in X^{**}$} with $\xi \cdot a = a\cdot \xi$, {for every $a$ in $A$}, and $d: A\to X$ is a linear map satisfying $d(b)\cdot a + b\cdot d(a)=  - 2 \xi \cdot (b a),$ for all $a,b\in A$ with $a b =0$. In such a case, given $a,b\in A$ with $a b =0,$ it is easy to see that $$ T(b) \cdot a + b\cdot T(a) = d(b)\cdot a + b\cdot d(a) + 2 \xi \cdot (b a) = 0.$$ Note that no other property on $d$ is employed. \smallskip

Suppose now that every Jordan derivation from $A$ to $X$ is a derivation. \smallskip

$(ii) \Rightarrow (iii)$ By hypotheses, we can write $T(a) = d(a) + \xi\cdot a$, for all $a\in A$, where {$\xi\in X^{**}$} satisfies $\xi\cdot a = a\cdot \xi$ for all $a\in A$, and $d: A\to X$ is a derivation satisfying $d(b)\cdot a + b \cdot d(a) = -2 \xi \cdot (b a)$ for all $a,b\in A$ with $a b = 0$. It is known that under these conditions $T$ is a generalized derivation, or alternatively, it can be directly checked that the identity
 \begin{equation}\label{equa T1}
	T(a b) = T(a) \cdot b +a \cdot T(b) - a \cdot \xi \cdot b,
\end{equation} holds for all $a,b \in A.$ Since \eqref{eq T is a kind of generalized antider} also holds, we arrive to $$T([a,b]) = - \xi \cdot [a,b], \hbox{ and consequently, } d([a,b])= -2 \xi \cdot [a,b], \hbox{ for all } a,b\in A.$$ 

$(iii) \Rightarrow (i)$ Fix any $a,b \in A$ with $a b =0$. Since $d([a,b]) = -2 \xi \cdot [a,b]$, it follows that $d(b a) = -2 \xi \cdot (b a)$. Therefore,
    \begin{align*}
        T(b)\cdot a + b\cdot T(a) &= (d(b)+ \xi \cdot b)\cdot a +b \cdot (d(a)+ \xi \cdot a)  \\
        &= d(b) \cdot a + b \cdot d(a) + 2 \xi \cdot (b a)\\
        &= d(b a) + 2 \xi \cdot (b a ) = 0,
    \end{align*} and hence $T$ is anti-derivable at zero.
\end{proof}

\begin{remark}\label{r the unital case in Thm 3.1 is trivial} We note that the difficulty in proving \Cref{thm anti-derivable at zero} arises from the absence of a unit element. Actually, if $A$ is an associative unital algebra satisfying property $\mathbb{B}$ and $X$ is a unital $A$-bimodule, the arguments given above, or even a simpler version of them, prove that the following statements are equivalent for every linear map $T: A \rightarrow X$:
\begin{enumerate}[$(i)$]
	\item $T$ is anti-derivable at zero.
	\item The element $T(\mathbf{1})\in X$ satisfies $T(\mathbf{1}) \cdot a = a \cdot T(\mathbf{1})$ for all $a\in A$, and  
	there exists a linear map (actually a Jordan derivation) $d: A\to X$ such that  
	$T(a) = d(a) +T(\mathbf{1}) \cdot a$, and $d(b)\cdot a + b\cdot d(a)=  - 2 T(\mathbf{1}) \cdot (b a),$ for all $a,b\in A$ with $a b =0$. 
\end{enumerate}
If we additionally assume that every Jordan derivation from $A$ to $X$ is a derivation, statements $(i)$ and $(ii)$ above are equivalent to:\begin{enumerate}[$(iii)$]
	\item The element $T(\mathbf{1})\in X$ satisfies $T(\mathbf{1}) \cdot a = a \cdot T(\mathbf{1})$ for all $a\in A$, and there exists a derivation $d: A\to X$ such that $T(a) = d(a) +T(\mathbf{1}) \cdot a$, and $d([a,b]) = -2T(\mathbf{1})\cdot [a,b]$, for all $a,b\in A$. In particular, $T$ is a generalized derivation.
\end{enumerate}	
\end{remark}

Anti-derivations and continuous linear maps which are anti-derivable at zero could simply reduce to the zero map in certain cases. For example, Theorem 4 in \cite{abulhamil2020linear} assures that for each C$^*$-algebra $A$, every anti-derivation $\delta : A\to X$ with $X = A, A^{*}$ or $A^{**}$ is zero. Furthermore, consider a unital Banach algebra $A,$ with unit $\mathbf{1},$ which is generated by commutators (see, for instance, the study in \cite{EusebioThiel2025rings}). Suppose additionally that every Jordan derivation from $A$ to a Banach $A$-bimodule $X$ is a derivation. Let $T: A\to X$ be a bounded liner map which is anti-derivable at zero. \Cref{thm anti-derivable at zero} implies that $T(a) = - \xi \cdot a$ and $d(a) = -2\xi \cdot a,$ for all $a \in A$, where $\xi\in X^{**}$ satisfies $\xi\cdot a = a\cdot \xi\in X$ for all $a\in A$. In particular, $0=d(\mathbf{1}) = -2 \xi \cdot \mathbf{1} = -2 \xi,$ since $d$ is a derivation, and thus $d(a) = 0 $ and $T(a)=0,$ for all $a \in A$. See \cite{Liu2025characterization} for additional examples of generalized matrix algebras on which every bounded linear map anti-derivable at zero is null. \smallskip

It should be highlighted that there exists a wide list of Banach algebras and modules satisfying that every Jordan derivation between them is a derivation. For example, every Jordan derivation from a C$^*$-algebra $A$ into a Banach $A$-bimodule is a derivation (see the introduction), and every Jordan derivation on a {CSL} algebra is a derivation (cf. \cite{Lu2009JordanCSLalgebras}). It is time to recall some notions about non-selfadjoint Banach algebras. Let $H$ be a separable complex Hibert space and $\mathcal{L}$ be a collection of closed subspaces of $H$. We say that $\mathcal{L}$ is a \emph{subspace lattice} on $H$ if $\{0\}$ and $H$ are both inside $\mathcal{L}$, and for every family $\{L_{i}\} \subseteq \mathcal{L}$ the intersection $\bigcap_{i} L_{i}$ and the closed linear span $\bigvee_{i} L_{i}$ belong to $\mathcal{L}$. We write $P_{L}$ for the orthogonal projection onto the subspace $L$. $\mathcal{L}$ is said to be a \emph{commutative subspace lattice} (\emph{CSL}, for short) if the projections in $\{P_{L}: L \in \mathcal{L}\}$ pairwise commute (i.e. $P_{L_{1}} P_{L_{2}} = P_{L_{2}} P_{L_{1}}$ for all $P_{L_{1}},P_{L_{2}} \in \mathcal{L}$). The subspace lattice algebra $Alg(\mathcal{L})$ corresponding to a subspace lattice $\mathcal{L}$ is defined by 
$$Alg(\mathcal{L}) \;:=\; \{\,T\in B(H): T(L)\subseteq L \text{ for all } L\in\mathcal{L}\,\},$$
that is, $T \in Alg(\mathcal{L})$ if, and only if, for every $L \in \mathcal{L}$, we have $TP_{L} = P_{L}T P_{L}$. Algebras of the form $Alg(\mathcal{L})$ are called reflexive operator algebras. A \emph{CSL algebra} is a reflexive operator algebra $A= Alg(\mathcal{L})$ whose lattice of invariant projections $\mathcal{L} = Lat (A)$ is a set of commuting projections. Finally, $A$ is called a \emph{CDCSL algebra} if it is a CSL algebra for which the lattice $\mathcal{L}$ is completely distributive as a lattice. It is known that all nest algebras and all complete atomic Boolean lattices are CDCSL algebras. 
We recommend \cite{LauLon1983rankone,Davidson1988nestalgs,ArgLamLon1991ABASL} and the references therein for more information about this topic.\smallskip 

It is worth pointing out some non-trivial examples of continuous linear maps which are anti-derivable at zero (see also \Cref{ex nontrivial anti-derivations,ex nontrivial anti-derivable at zero}). Every C$^*$-algebra satisfies all the hypotheses in \Cref{thm anti-derivable at zero} (cf. \cite[Theorem 3.2.21]{Dales2001banachalgs}, \cite[\S 1.3]{ABEV2009mapspreservingzpd}, and \cite{Johnson1996}), and so the result can be applied to this case. In \Cref{sec C-star-algebras} below we shall explore anti-derivable maps at zero on C$^*$-algebras and we shall obtain non-trivial examples even without assuming $X$ is unital. {For the moment, let us find some other applications of our result. In the following $Alg(\mathcal{L})$ denotes a CDCSL algebra.} Theorem 2 in \cite{LiPanSu2024CDCSLzpd} shows that $Alg(\mathcal{L})$ satisfies property $\mathbb{B}$. Moreover, every Jordan derivation from $Alg(\mathcal{L})$ into itself is a derivation (cf. \cite{Lu2009JordanCSLalgebras}). So, by considering a CDCSL algebra (or in particular, a nest algebra) $A$ and $X=A$, we are in a position to apply {\Cref{thm anti-derivable at zero}. 

\begin{corollary}\label{coro CDCSL antiderivable at zero}
    Let $A$ be a CDCSL algebra, and let $T: A \to A$ be a continuous linear map on $A$. Then the following statements are equivalent:
    \begin{enumerate}[$(i)$]
        \item $T$ is anti-derivable at zero.
        \item There exists an element $\xi $ in  $\mathcal{Z}(A)$ and a (continuous) derivation $d: A \rightarrow A$ satisfying $T(a) = d(a) +\xi  a$, $T([a,b]) = - \xi  [a,b]$ and $d([a,b]) = -2 \xi [a,b]$ for all $a,b \in A$. In particular, $T$ is a generalized derivation.
    \end{enumerate}
\end{corollary}

The next lemma {has been borrowed} from \cite[Lemma 3]{abulhamil2020linear}.

\begin{lemma}\label{AJP equi of anti-derivations}
    Let $d: A \rightarrow X$ be a linear map from an associative algebra into an $A$-bimodule. Then the following statements are equivalent:
    \begin{enumerate}[$(i)$]
        \item $d$ is a derivation and $d([a,b]) = 0,$ for all $a,b \in A$;
        \item $d$ is an anti-derivation and $d([a,b]) = 0,$ for all $a,b \in A$.
    \end{enumerate}
\end{lemma}

In view of \Cref{ex nontrivial anti-derivable at zero}, we cannot always expect that the derivation $d$ appearing  in \Cref{thm anti-derivable at zero}$(iii)$ is an anti-derivation. We characterize next when this can happen. 

\begin{corollary}\label{coro equiva of d becomes an anti-derivation}
Under all the assumptions stated at \Cref{thm anti-derivable at zero}, the derivation $d$ in statement $(iii)$ is an anti-derivation if, and only if, $\xi \cdot [a,b] = 0,$ for all $a,b \in A$ if, and only if, $d([a,b]) = 0,$ for all $a,b \in A$.
\end{corollary}

\begin{proof} Suppose first that $d$ is an anti-derivation, and thus $$d(a b) = d(a) \cdot b + a \cdot d(b) = d(b a) ,$$  for all $a,b \in A$. This implies that $0= d([a,b]) = -2 \xi \cdot [a,b]$, equivalently, $\xi \cdot [a,b] = 0,$ for all $a,b \in A$. Conversely, if $\xi \cdot [a,b] = 0,$ for all $a,b \in A$, then $d([a,b]) = -2 \xi \cdot [a,b] = 0$. Clearly $d$ is an anti-derivation by \Cref{AJP equi of anti-derivations}.\smallskip 

Finally, $d([a,b]) = 0,$ for all $a,b \in A$ if, and only if, $\xi \cdot [a,b] = 0,$ for all $a,b \in A$ since $d([a,b]) = -2 \xi \cdot [a,b]$.
\end{proof}

Let $X$ be an $A$-bimodule over an associative algebra $A$.  Recall that a derivation $d:A\to X$ is called {\emph{inner}} if there exists $x_0\in X$ satisfying $d(a) =[a,x_0] = a x_0 - x_0 a,$ for all $a\in A$.\smallskip 

If in \Cref{thm anti-derivable at zero} we assume that every Jordan derivation from $A$ to $X$ is an inner derivation the conclusion is a bit stronger. 

\begin{proposition}\label{prop suppose inner der charac}
Let $A$ be a Banach algebra with a bounded approximate unit and  satisfying property $\mathbb{B}$, and let $X$ be an essential Banach $A$-bimodule such that every Jordan derivation from $A$ to $X$ is an inner derivation. Suppose that $T: A\to X$ is a continuous linear map. Then $T$ is anti-derivable at zero if, and only if, there exist $u\in X$ and $v \in X^{**}$ such that $(u-v) \cdot a = a\cdot (u-v)\in X$, $T(a) = a\cdot u - v\cdot a$, and $[b a, u ] + 2 (u-v) \cdot (b a) =0,$ for all $a,b\in A$ with $a b =0$. 
\end{proposition}

\begin{proof} Suppose first that $T$ is anti-derivable at zero. It follows from \Cref{thm anti-derivable at zero} that there exist a continuous derivation $d: A\to X$ and an element $\xi \in X^{**}$ such that $\xi \cdot a = a\cdot \xi\in X,$ $T(a) = d(a) + \xi\cdot a,$ and $T([a,b]) = -\xi \cdot [a,b]$, for all $a,b \in A$.\smallskip
	
Since, by hypotheses, $d$ is an inner derivation, there exists $u \in X$ such that $d(a) =[a,u] = a u - u a,$ for all $a \in A$. By defining $v = u - \xi \in X^{**}$ we have $T(a) = a\cdot u - v\cdot a,$ and $(u - v)\cdot a = a\cdot (u-v)\in X$, for all $a \in A$.\smallskip

Take now $a,b\in A$ with $a b =0$. The identity $d([a,b]) = -2 \xi\cdot[a,b] = -2 [a,b]\cdot \xi$, assures that $[- b a, u] = 2 (ba)\cdot (u-v)$, which proves the desired properties for $u$ and $v$.\smallskip   

Conversely, if there exist $u\in X$ and $v \in X^{**}$ such that $(u-v) \cdot a = a\cdot (u-v)\in X$, $T(a) = a\cdot u - v\cdot a$, and $[b a, u ] + 2 (u-v) \cdot (b a) =0,$ for all $a,b\in A$ with $a b =0$, it is easy to see that 
  $$\begin{aligned}
  	T(b)&\cdot a + b\cdot T(a) = (b\cdot u - v\cdot b)\cdot a  + b\cdot (a\cdot u - v\cdot a)\\
  	& = b\cdot u \cdot a -  v\cdot (b a) + (b a)\cdot u  - b\cdot v\cdot a  \\
  	&= b\cdot u \cdot a -  u\cdot (b a) +  (u-v) \cdot (b a) + (b a)\cdot u  - b \cdot u\cdot a + b \cdot (u-v) \cdot a \\
  	&= [b a , u] + 2 (u-v) \cdot (b a) =0, 
  \end{aligned}$$ as desired.   
\end{proof}

Recall that every von Neumann algebra $A$ is unital, every derivation on $A$ is inner \cite[Theorem 4.1.6]{SakaiBook71}, and every Jordan derivation on $A$ is a derivation. So, the previous proposition applies when $A= X$ is a von Neumann algebra. In \Cref{sec C-star-algebras} below we shall improve this conclusion.\smallskip

Since every derivation on nest algebra is inner (see for instance, \cite[Lemma 2.4]{HanWei1995local}), and each nest algebra is a CDCSL algebra, it follows from \cite{Lu2009JordanCSLalgebras} that
every Jordan derivation on a nest algebra is inner. Moreover every nest algebra satisfies property $\mathbb{B}$ \cite[Theorem 2]{LiPanSu2024CDCSLzpd}.

\begin{corollary} Let $A$ be a nest algebra, and let $T$ be a continuous linear map on $A$. Then $T$ is anti-derivable at zero if, and only if, there exist $u,v \in A$, such that $u-v \in \mathcal{Z}(A)$, $T(a) = a u - v a$, and $[b a, u ] + 2 (u-v) b a  =0,$ for all $a,b\in A$ with $a b =0$.
\end{corollary}

Now let $A$ be an arbitrary unital {algebra}. It follows from \cite[Theorem 2.1]{HuangLiQian2020derivations} that every derivation $D$ on the matrix algebra $M_{n}(A)$ ($n \geq 2$) can be written in the form $D = D_{a} + \Tilde{d},$ where $D_{a}$ is an inner derivation associated with some element $a \in M_n(A),$ $d :A \to A$ is a derivation, and $\Tilde{d}$ is the derivation on $M_{n}(A)$ induced by $d$ by the assignment $\Tilde{d} \left( (a_{ij} ) \right) =(d(a_{ij}))$. It is known that $M_{n}(A)$ satisfies property $\mathbb{B}$ for all $n\geq 2$ (see \cite[Theorem 2.1]{Bresar2009zpdmatrixalgebras}). If we additionally assume that every derivation on $A$ is inner, then $\Tilde{d}$ must be inner and hence $D$ is inner too. 
When in the proof of \Cref{prop suppose inner der charac} we replace \Cref{thm anti-derivable at zero} with \Cref{r the unital case in Thm 3.1 is trivial} we arrive to the next result.

\begin{corollary} Let $A$ be a complex algebra with unit $\mathbf{1}$, and let $T$ be a linear map on the matrix algebra $M_{n}(A)$. Suppose that every derivation on $A$ is inner.
Then $T$ is anti-derivable at zero if, and only if, there exist $u,v \in M_{n}(A)$, such that $u-v \in \mathcal{Z}(M_{n}(A))$, $T(a) = a u -v a$, and $[b a, u ] + 2 (u-v) b a =0,$ for all $a,b \in A$ with $a b =0$.
\end{corollary}

\section{Linear maps anti-derivable at zero whose domain is a C$^*$-algebra}\label{sec C-star-algebras}

Throughout this section, let $X$ be an essential Banach $A$-bimodule over a C$^*$-algebra $A$. It is well known that we can always find a bounded approximate unit $(e_j)_j$ in $A$ which converges in the weak$^*$-topology of $A^{**}$ to the unit element $\mathbf{1}\in A^{**}$ \cite[Theorem 3.2.21]{Dales2001banachalgs}. In this case we have $\mathbf{1} \cdot \xi = \xi \cdot \mathbf{1} = \xi$ for all $\xi \in X,$ and $\eta \cdot \mathbf{1} = \eta$ for all $\eta \in X^{**}$ (cf. \Cref{r essential modules and bau} and comments before \cite[Lemma 5]{abulhamil2020linear}).\smallskip

Our study begins with a technical observation. 

\begin{lemma}\label{lem eta commutes with A**} Let $A$ be a Banach algebra admitting a bounded approximate unit $(e_j)_j$, and let $X$ be an essential Banach $A$-bimodule. Suppose $T:A \to X$ is a continuous generalized derivation, i.e., a bounded linear map for which there exists $\xi \in X^{**}$ satisfying \begin{equation}\label{lem equa generalized der in eta x x eta}
        T(a b) = T(a)\cdot b + a \cdot T(b) - a \cdot \xi \cdot b,
\end{equation} 
for all $a,b \in A$. Suppose further {that} $\xi \cdot a = a \cdot \xi$ for every $a\in A$, and let $e$ be a projection obtained as a weak$^*$ cluster point of $(e_j)_j$ in $A^{**}$. Then the element $\eta := e \cdot \xi \cdot e= e \cdot \xi  \in X^{**}$ satisfies $\eta \cdot a = a \cdot \eta$, for all $a \in e A^{**}.$ If $A$ is Arens regular the conclusion actually holds for all $a\in A^{**}$. 
\end{lemma}

\begin{proof} By replacing $b$ with $e_j$ in \eqref{lem equa generalized der in eta x x eta} and applying the hypotheses on $\xi$ we get $$T(a e_j ) = T(a)\cdot e_j + a \cdot T(e_j) - a \cdot \xi \cdot e_j = T(a)\cdot e_j + a \cdot T(e_j) - (a  e_j) \cdot \xi,$$ for all $j$. So, taking weak$^*$ and norm limits of an appropriate subnet, and having in mind that the map $A^{**}\mapsto A^{**}$, {$z\mapsto a z$} is weak$^*$-continuous for all $a\in A$, we deduce that 
    \begin{equation*}
        T(a) = T(a) + a \cdot T^{**}(e) - a \cdot \xi,
    \end{equation*} equivalently, 
    \begin{equation}\label{lem equa x xi = x T1}
        a \cdot T^{**}(e) = a \cdot \xi, \hbox{ for every } a\in A.
    \end{equation}

Next, taking $a = e_j$ in the previous identity and employing the same technique above we get $e \cdot \xi  = e \cdot T^{**}(e)$. So $$\eta=e \cdot \xi \cdot e = e \cdot T^{**}(e) \cdot e.$$ 
    
By Goldstine's theorem, for each $a\in A^{**}$, there is a net $(a_{\mu})_{\mu}$ in $A$ which converges in the weak$^*$-topology of $A^{**}$ to $a$. Since
    $$T(e_j a_{\mu}) = T(e_j)\cdot a_{\mu} + e_j \cdot T(a_{\mu}) - e_j \cdot \xi \cdot a_{\mu},$$ we have 
    \begin{equation}\label{eq 410}
    e_j \cdot \xi \cdot a_{\mu} = T(e_j)\cdot a_{\mu} + e_j \cdot T(a_{\mu}) - T(e_j a_{\mu}).	
    \end{equation} Note that the properties of the first Arens product lead to 
    $$w^{*}\text{-}\lim_{j} \lim_{\mu} T(e_j)\cdot a_{\mu} = T^{**}(e) \cdot a, \ \ w^{*}\text{-}\lim_{j} \lim_{\mu} e_j \cdot T(a_{\mu}) = e \cdot T^{**}(a)$$ and
    $$w^{*}\text{-}\lim_{j} \lim_{\mu}T(e_j a_{\mu}) = T^{**}({\color{red}e} a).$$
    Furthermore, by a new application of the first Arens product's properties and the observation before \Cref{r essential modules and bau}, we obtain 
    \begin{align*}
        w^{*}\text{-}\lim_{j} \lim_{\mu} e_j \cdot \xi \cdot a_{\mu} &= w^{*}\text{-}\lim_{j} \lim_{\mu} (e_j  a_{\mu}) \cdot \xi = w^{*}\text{-}\lim_{j} (e_j a)\cdot \xi\\
        &=\left(w^{*}\text{-}\lim_{j} e_j a\right) \cdot \xi = (e a) \cdot \xi.
    \end{align*}
Therefore, the identity in \eqref{eq 410} assures that 
 \begin{equation}\label{eq 0410 b} {(e a)} \cdot \xi =  T^{**}(e) \cdot a +  e \cdot T^{**}(a) -  T^{**}( {\color{red}e} a ) \ \ \ (\forall a\in A^{**}).
 \end{equation}
    
Suppose, finally, that $a\in e A^{**}$. By multiplying the identity in \eqref{eq 0410 b} on both sides by $e$ we obtain 
$$\begin{aligned}  a \cdot \eta = e \cdot ((a e)  \cdot \xi) \cdot e &=
e \cdot (a \cdot \xi) \cdot e\\
&= e\cdot (T^{**}(e) \cdot a) \cdot e +  e\cdot (e \cdot T^{**}(a))\cdot e - e\cdot  T^{**}(a) \cdot e \\
&= (e\cdot T^{**}(e) \cdot e) \cdot a = (e\cdot \xi \cdot e) \cdot a = \eta \cdot a,
\end{aligned}$$
which completes the proof.
\end{proof}

As we commented in the introduction, bounded linear operators from a C$^*$-algebra $A$ to a Banach $A$-bimodule which are anti-derivable at zero were completely determined in \cite[Theorem 6]{abulhamil2020linear}. Unfortunately, one of the steps in the arguments relies on a property which is not, in general, true. Namely, the claim that the opposite module, $X^{op},$ of a Banach $A$-bimodule $X$ is a Banach $A$-bimodule, may fail in some cases. Recall that the module products on $X^{op}$ are defined by $a\odot x = x\cdot a$ and $x\odot a = a\cdot x$, for all $a\in A$, $x\in X$. The identity $a \odot ( b \odot x) = (ab) \odot x$ does not necessarily hold when $A$ is not commutative. This difficulty affects the proof of the just quoted result. In our next theorem we show that the original statement in \cite[Theorem 6]{abulhamil2020linear} is true by providing a complete new proof of the result, which requires a more elaborated argument.  

\begin{theorem}\label{thm C-star algebra being anti-derivable at zero} Let $T: A \rightarrow X$ be a continuous linear operator, where $A$ is a C$^*$-algebra and $X$ is an essential Banach $A$-bimodule. Then the following are equivalent:
\begin{enumerate}[$(i)$]
\item $T$ is anti-derivable at zero.
\item There exist an element $\eta \in  X^{**}$ and an anti-derivation $d: A \rightarrow X$ satisfying $\eta \cdot a = a \cdot \eta \in X$, $\eta \cdot [a,b] = 0$ {\rm(}i.e., $L_{\eta}: A \to A$, $L_{\eta} (a) = \eta \cdot a$ vanishes on commutators{\rm)}, and $T(a) = d(a) +\eta \cdot a$, for all $a,b \in A$. 
\end{enumerate}
Furthermore, in case that $T$ is anti-derivable at zero, the element $\eta$ in statement $(ii)$ actually satisfies that $\eta = z_{I_1} \cdot \eta = \eta \cdot z_{I_1}$, where $z_{I_1}$ is the central projection in $A^{**}$ satisfying that $z_{I_1} A^{**}$ is the type $I_1$ part of $A^{**}$. If $A$ is unital $T(\mathbf{1})= \eta\in X$. 
\end{theorem}

\begin{proof}
$(ii) \Rightarrow (i)$ Fix $a,b \in A$ with $a b =0$. Since $d$ is an anti-derivation, $\eta$ commutes with all elements in $A$, and $0= \eta \cdot [a,b] = - \eta \cdot (b a) $, it follows straightforwardly that  
    \begin{align*}
        0=T(a b) &= d(a b) +\eta \cdot (ab)  = d(a b)= d(b)\cdot a + b\cdot d(a)  \\
        &= T(b) \cdot a -\eta \cdot (ba) + b\cdot T(a) -\eta \cdot (ba)  \\
		&=T(b)\cdot a + b\cdot T(a).
    \end{align*}
    So $T$ is anti-derivable at zero.\smallskip

$(i) \Rightarrow (ii)$ We can clearly apply \Cref{thm anti-derivable at zero} $(i)\Leftrightarrow (iii)$ to this special case. Therefore there exist an element $\xi \in  X^{**}$ and a derivation $d: A \rightarrow X$ satisfying $T(a) = d(a) +\xi \cdot a$, $\xi \cdot a = a \cdot \xi \in X$,  $T([a,b]) = - \xi \cdot [a,b]$ and $d([a,b]) = -2 \xi \cdot [a,b],$ for all $a,b \in A$. Consequently, \begin{equation}\label{thm equa in C star alg T generalized der}
	T(a b) = T(a)\cdot b + a \cdot T(b) - a \cdot \xi \cdot b,
\end{equation} for all $a,b \in A$. Let $\mathbf{1}$ denote the unit of $A^{**}$. Since every C$^*$-algebra is Arens regular and admits a bounded approximate unit, \Cref{lem eta commutes with A**} implies that $\eta = \mathbf{1}\cdot \xi \cdot \mathbf{1}$ commutes with all elements in $A^{**}$. It follows from \eqref{thm equa in C star alg T generalized der} that $$\begin{aligned}
T(a b) &= T(a)\cdot b + a \cdot T(b) - (a \mathbf{1})\cdot \xi \cdot (\mathbf{1} b) \\
&= T(a)\cdot b + a \cdot T(b) - a \cdot \eta \cdot b \\
&= T(a)\cdot b + a \cdot T(b) - (a b) \cdot \eta
\end{aligned}$$ for all $a,b\in A$.  Now, by combining the weak$^*$-density of $A$ in $A^{**}$, the weak$^*$-continuity properties of $T^{**}: A^{**}\to X^{**}$ and of the $A^{**}$-module operations of $X^{**}$, we can easily obtain that the identity \begin{equation}\label{thm equa in C star alg T** generalized der}
T^{**}(a b) = T^{**}(a) \cdot b + a \cdot T^{**}(b) - (a b) \cdot \eta =  T^{**}(a) \cdot b + a \cdot T^{**}(b) - a \cdot \eta \cdot b 
\end{equation}
holds for all $a,b \in A^{**}$. We can also arrive to the previous identity by just applying that $a\cdot \xi\in X$ for all $a\in A$, and a similar approach via weak$^*$-limits.\smallskip
    
The bi-transpose of $d$, $d^{**}: A^{**} \rightarrow X^{**},$ is a (continuous) derivation on $A^{**}$ and satisfies $d^{**}(a) = T^{**}(a) - a\cdot \eta = T^{**}(a) - \eta \cdot a $. By combining the Arens regularity of $A$ with the weak$^*$-continuity properties of $T^{**}$ and of the $A^{**}$-bimodule operations on $X^{**},$ and with the identities $$T([a,b]) = - \xi \cdot [a,b]= - [a,b] \cdot \xi = - a \xi b + b  \xi a = - a \eta b + b  \eta a  = - [a,b] \cdot \eta,$$ and $$d([a,b]) = -2 [a,b]  \cdot \xi = -2 [a,b]  \cdot \eta,$$ for all $a,b \in A$, we arrive to  
\begin{equation}\label{eq T and d bitransposed on commutators} T^{**}([a,b]) = - [a,b] \cdot \eta, \text{ and } d^{**}([a,b]) = -2 [a,b] \cdot \eta,
\end{equation} for all $a,b \in A^{**}$. \smallskip 

By \Cref{AJP equi of anti-derivations}, in order to prove that $d$ is an anti-derivation it suffices to show that $d^{**}([a,b]) = -2 [a,b] \cdot \eta = 0,$ for every $a,b \in A^{**}$ (which actually shows the stronger conclusion that $d^{**}$ is an anti-derivation). \smallskip

The structure theory of von Neumann algebras (cf. \cite[\S V]{Tak79}) assures that $A^{**}$ is uniquely decomposable into 
a direct sum of of the form 
$$A^{**} = p_{1}A^{**} \mathop{\bigoplus}\limits^{\ell_{\infty}} p_{2}A^{**} \mathop{\bigoplus}\limits^{\ell_{\infty}} p_{3}A^{**},$$
where $p_{1}, p_{2}, p_{3}$ are pairwise orthogonal central projections in $A^{**}$ such that $p_{1}+ p_{2}+ p_{3} =\mathbf{1}$ and $p_{1} A^{**}$ is of type I finite, $p_{2}A^{**}$ is of type II$_{1},$ and $p_{3} A^{**}$ is properly infinite.\smallskip

For each $i\in \{1,2,3\},$ define $d^{**}_{i} : p_{i} A^{**} \rightarrow p_{i} \cdot X^{**}$ by $d_i (a) = p_{i} \cdot d^{**}(a)$. It can be easily checked that $d^{**}_{i}$ is a derivation.\smallskip

It follows from a celebrated result by Fack and de la Harpe in \cite[Theorem 3.2 and Theorem 3.10]{FackHarpe1980sommes} that every element with zero trace in a finite von Neumann algebra $W$ can be written as the sum of ten commutators in $W$, and moreover, every element in a properly infinite von Neumann algebra can be expressed as a sum of two commutators. So, each $a\in p_{3} A^{**}$ can be written as the sum of two commutators, and thus $d^{**}_{3} (a) = -2 a \cdot \eta$ for all $a\in p_{3} A^{**}$ (cf. \eqref{eq T and d bitransposed on commutators}). Since $p_3$ is the unit element of $p_3 A^{**}$ and $p_3$ is a left unit for $p_3\cdot X^{**}$ we must have   
    $$\begin{aligned}
    	d^{**}_{3}(p_{3}) &= d^{**}_{3}(p_{3})\cdot p_3 + p_3 \cdot d^{**}_{3}(p_{3})\\
    	&= -2 p_{3} \cdot \eta \cdot p_3 + d^{**}_{3}(p_{3}) = -2 p_{3} \cdot \eta + d^{**}_{3}(p_{3}),
    	\end{aligned}$$ which assures that  \begin{equation}\label{equa dstar_3 xi p_3}  p_{3} \cdot \eta =  \eta \cdot p_{3}= 0.
    \end{equation}
    
Let $\tau: p_{1}A^{**} \mathop{\bigoplus}\limits^{\ell_{\infty}} p_{2}A^{**} \to \mathcal{Z}\left(p_{1}A^{**} \mathop{\bigoplus}\limits^{\ell_{\infty}} p_{2}A^{**}\right)$ denote the faithful center-valued trace on the finite von Neumann algebra $p_{1}A^{**} \mathop{\bigoplus}\limits^{\ell_{\infty}} p_{2}A^{**}$.\smallskip 

On the type II$_1$ von Neumann algebra $p_2 A^{**}$, we can apply the Halving lemma (see \cite[Proposition V.1.35]{Tak79}) to deduce the existence of two orthogonal projections $q_{1},q_{2} \in p_{2} A^{**}$ which are (Murry-von Neumann) equivalent and satisfy $p_{2} = q_{1} +q_{2}$ and $\tau(q_{1}) = \tau (q_{2})$. Hence $\tau (q_{1} -q_{2}) = 0,$ and so $q_{1} -q_{2}$ writes as a finite sum of commutators in $p_{2} A^{**}$ (cf. \cite[Theorem 3.2]{FackHarpe1980sommes}). Consequently, by \eqref{eq T and d bitransposed on commutators}, $d^{**}_{2}(q_{1}-q_{2}) = -2 (q_{1}-q_{2}) \cdot \eta.$ We therefore have $d^{**}_{2}(p_{2}) = d^{**}_{2}(p_{2})\cdot p_2 + p_2\cdot d^{**}_{2}(p_{2}) = d^{**}_{2}(p_{2})\cdot p_2 + d^{**}_{2}(p_{2})$  and  
    \begin{align*}
        0 = d^{**}_{2}(p_{2})\cdot p_2 &= d^{**}_{2}((q_{1}-q_{2})^{2})\cdot p_2\\
        &= d^{**}_{2}(q_{1}-q_{2}) \cdot (q_{1}-q_{2}) + (q_{1}-q_{2}) \cdot d^{**}_{2}(q_{1}-q_{2}) \cdot p_2 \\
        &= -4\eta \cdot (q_{1}-q_{2})^{2} = -4\eta \cdot p_{2}.
    \end{align*}
We have shown that  \begin{equation}\label{equa dstar_2 xi p_2}
        0=  \eta \cdot p_{2} = p_{2} \cdot \eta.
    \end{equation}
    
Finally, we analyse the structure of the type I finite von Neumann algebra $p_1 A^{**}$. It is known that there exists a sequence $\{z_{j}\}_{j}$ of pairwise orthogonal central projections in $p_{1} A^{**}$ such that 
    $\displaystyle p_{1}A^{**} \cong \mathop{\bigoplus}\limits^{\ell_{\infty}}_{j\in N_0} z_{j}A^{**}$  with $\displaystyle \sum_{j}z_{j} = p_{1}$ and $N_0\subseteq \mathbb{N}$, where for each $j$, $z_{j}A^{**} $ is a type $I_{n_j}$ von Neumann algebra, and thus $^*$-isomorphic to $C(K_j, B(\ell_2^{n_{j}})),$ for some hyperstonean space $K_j$ and some $n_{j} \in \mathbb{N}$ with $n_{j_1}\neq n_{j_2}$ for $j_1\neq j_2$ (cf. \cite[Theorem V.1.27]{Tak79}).\smallskip

Fix an arbitrary $j\in N_0$ and an orthonormal basis $\{\xi_{1},\cdots, \xi_{n_{j}} \}$ of the Hilbert space $\ell_2^{n_{j}}$. Given $\xi,\eta\in \ell_2^{n_{j}}$, we shall write $\eta \otimes \xi$ for the operator in $B(\ell_2^{n_{j}})$ defined by $\eta \otimes \xi (\zeta) = \langle \zeta| \xi\rangle \eta.$ For each $ 1\leq i < k \leq n_{j}$ we shall write $u_{ik} $ for the element in $z_{j}A^{**} \cong C(K_j, B(\ell_2^{n_{j}}))$ defined as the constant function with constant value $\xi_{i} \otimes \xi_{k} +\xi_{k} \otimes \xi_{i}$. \smallskip

As in the arguments in the previous paragraphs, the map $D_{j}:= z_j\cdot d_{1}^{**}\arrowvert_{z_{j}A^{**}}: z_{j} A^{**} \to z_{j} \cdot X^{**}$ is a derivation. It follows that $D_j (z_j) = z_j \cdot D_j (z_j) + D_j (z_j) \cdot z_j  = D_j (z_j) + D_j (z_j) \cdot z_j,$ and thus $D_j (z_j) \cdot z_j =0$.\smallskip

If $n_{j}$ is even, the element $u_{2\ell-1,2\ell}$ has zero trace (i.e. $\tau(u_{2\ell-1,2\ell}) = 0$) for all $1 \leq \ell \leq \frac{n_{j}}{2}$, then it follows from \cite[Theorem 3.2]{FackHarpe1980sommes} that $u_{2\ell-1,2\ell}$ writes as a finite sum of commutators in $z_j A^{**}$, and thus by \eqref{eq T and d bitransposed on commutators} we get $$D_{j}(u_{2\ell-1,2\ell}) = -2\eta \cdot u_{2\ell-1,2\ell},$$ and consequently, $D_{j}(u_{2\ell-1,2\ell}^{2}) = -4\eta \cdot u_{2\ell-1,2\ell}^{2},$ since $D_{j}$ is a derivation and $\eta$ commutes with every element in $A^{**}$.  
Since $D_j (z_j) \cdot z_j =0$, we deduce that 
    \begin{align*}
        0= D_j (z_j) \cdot z_j  &= D_{j}\left(\mathop{\sum}\limits_{1 \leq \ell \leq \frac{n_{j}}{2}} u_{2\ell-1,2\ell}^{2}\right)\cdot z_j = \mathop{\sum}\limits_{1 \leq \ell \leq \frac{n_{j}}{2}} D_{j}( u_{2\ell-1,2\ell}^{2}) \cdot z_j  \\ 
        &= -4 \eta \cdot \left(\mathop{\sum}\limits_{1 \leq \ell \leq \frac{n_{j}}{2}} u_{2\ell-1,2\ell}^{2} \right) \cdot z_j  = -4 \eta \cdot z_{j},
    \end{align*} which proves that $\eta \cdot z_j = z_j \cdot \eta$ whenever $n_j$ is even.\smallskip
    
If $n_{j}$ is odd with $n_j\geq 3$, as before, $\tau(u_{ik}) = 0$ for every $1\leq i < k \leq n_{j},$ and hence $D_{j}(u_{ik}) = -2 \eta \cdot u_{ik}$ by \cite[Theorem 3.2]{FackHarpe1980sommes} and \eqref{eq T and d bitransposed on commutators}. Furthermore, $D_{j}(u_{ik}^{2}) = -4\eta \cdot u_{ik}^{2}$ since $D_{j}$ is a derivation. A simple calculation shows that  
    $$\mathop{\sum}\limits_{1\leq i < k \leq n_{j}} u_{ik}^{2} = (n_{j} -1) z_{j}.$$  We therefore obtain 
    \begin{align*}
      0&= (n_{j}-1)D_{j}(z_{j})\cdot z_j =  D_{j}\left( \mathop{\sum}\limits_{1\leq i < k \leq n_{j}}u_{ik}^{2} \right) \cdot z_j \\
      &=  \left(\mathop{\sum}\limits_{1\leq i < k \leq n_{j}} D_{j}(u_{ik}^{2}) D_j (z_j)\right) \cdot z_j  = -4\eta \cdot \left(\mathop{\sum}\limits_{1\leq i < k \leq n_{j}} u_{ik}^{2}\right) \cdot z_j \\
      &= -4(n_{j}-1)\eta \cdot z_{j},
    \end{align*} which implies that $\eta \cdot z_{j} =  z_{j} \cdot  \eta =0$ if $n_{j}$ is odd with $n_j\geq 3$.\smallskip

In case that there exists $j_1\in N_0$ with $n_{j_1} = 1$, the von Neumann algebra $z_{j_1} A^{**} \cong C(K, \mathbb{C})$ is abelian, so we obviously have $d^{**}([x,y]) = 0,$ for all $x,y \in z_{j_1}A^{**}$. It can be easily concluded that  
   \begin{equation}\label{equa dstar_1 xi p_1}
       p_1 \cdot \eta= \eta \cdot p_{1} =\eta \cdot z_{j_{1}} +\sum_{n_{j} \text{ is even}} (\eta \cdot z_{j}) +\sum_{n_{j}\neq 1 \text{ is odd}} (\eta \cdot z_{j}) = \eta \cdot z_{j_{1}} = z_{j_{1}} \cdot \eta,\end{equation} and thanks to \eqref{equa dstar_1 xi p_1}, \eqref{equa dstar_2 xi p_2}, and \eqref{equa dstar_1 xi p_1}, we arrive to 
   \begin{align*}
       \eta = \eta \cdot \mathbf{1} = \eta \cdot (p_{1} +p_{2} +p_{3})=\eta \cdot z_{j_{1}}.  
   \end{align*} Finally, given $a,b \in A^{**},$ it easily follows from \eqref{eq T and d bitransposed on commutators}, the fact that $z_{j_1}$ is central, the summand $z_{j_1} A^{**}$ is abelian, and the previous identity that 
   \begin{align*}
       d^{**}([a,b]) &= -2 [a,b] \cdot \eta = - 2 \eta\cdot [a,b]  = -2(\eta \cdot z_{j_{1}})\cdot [a,b]\\
       &= -2\eta \cdot (z_{j_{1}}[a,b]) = -2\eta \cdot [z_{j_{1}}a, z_{j_{1}}b] = 0,
   \end{align*} which finishes the proof of $(i)\Rightarrow (ii)$.\smallskip
   
The first part of the final comments has been seen above, while for the second part we simply observe that $d(\mathbf{1}) =0$ and $T(\mathbf{1}) = \eta \cdot \mathbf{1} =\mathbf{1} \cdot \eta = \eta$. 
\end{proof}

The following result is an immediate consequence of our previous \Cref{thm C-star algebra being anti-derivable at zero}.

\begin{corollary}
Let $T: A \rightarrow X$ be a continuous linear operator, where $A$ is a C$^*$-algebra and $X$ is an essential Banach $A$-bimodule. Suppose $A^{**}$ contains no type $I_{1}$ summand. Then the following are equivalent:
    \begin{enumerate}[$(i)$]
        \item $T$ is anti-derivable at zero.
        \item $T$ is an anti-derivation.
    \end{enumerate}    
\end{corollary}

\begin{proof} Clearly, we only need to prove that $(i)\Rightarrow (ii)$. If we assume that $T$ is anti-derivable at zero. Let $\eta$ be the element give by \Cref{thm C-star algebra being anti-derivable at zero}$(ii)$. The just quoted result actually assures that $\eta = \eta \cdot z_{I_1} =z_{I_1}\cdot \eta$, where $z_{I_1}$ is the central projection in $A^{**}$ which determines the type $I_{1}$ part of this von Neumann algebra. So, our assumptions imply that $z_{I_1}=0$ and hence $\eta = 0$. Consequently, $d^{**}(a) = T^{**}(a) -\eta \cdot a =  T^{**}(a),$ for all $a\in A^{**},$ which shows that $T^{**}= d^{**}$ is an anti-derivation. In particular, $T : A \to X$ is an anti-derivation.
\end{proof}

The study of linear maps on von Neumann algebras which are anti-derivable at zero can be done without assuming their continuity. However, it should be noted that every anti-derivation on a von Neumann algebra is zero (cf. \cite[Theorem 4]{abulhamil2020linear}).

\begin{corollary} Let $T: M \to M$ be a linear map on a von Neumann algebra $M$. Then the following are equivalent:
\begin{enumerate}[$(i)$]
\item $T$ is anti-derivable at zero.
\item There exists an element $\eta$  in the type $I_1$ part of $M$ satisfying $T(a) = \eta  a$ for all $a\in M$. 
\end{enumerate} 
\end{corollary}

\begin{proof}
Suppose first that $T$ is anti-derivable at zero. Fix any $a,b,c \in M_{sa}$ such that $ab = bc =0$.  Then 
    $$T(a) b + a T(b) = 0,$$
since $ba = 0$ and $T$ is anti-derivable at zero. Hence
    $$a T(b)c = (-T(a)b)c = -T(a) (bc) = 0 .$$
It follows from \cite[Corollary 2.15 and Corollary 2.13]{ABAPe2018linearmaps} that $T$ is a continuous generalized derivation on $M$. It follows that $T$ is a bounded linear map which is anti-derivable at zero, thus we deduce from \Cref{thm C-star algebra being anti-derivable at zero} that there exist an element $\eta \in \mathcal{Z}(M)$ and an anti-derivation $d$ on $M$ such that $T(a) = d(a) + \eta a,$ and $\eta [a,b] = 0$
for all $a,b \in M$. Having in mind that every anti-derivation on $M$ is zero (cf. \cite[Theorem 4]{abulhamil2020linear}), we obtain 
$T(a) = \eta a,$ for all $a\in A$.\smallskip
    
We shall next show that $\eta$ lies in the type $I_1$ part of $M$. Observe that $\eta [a,b] =0$ for all $a,b\in M$. We employ, once again, the Murray-von Neumann classification to decompose $M$ in the form $M = p_1 M + p_2 M + p_3 M$, where $p_1,p_2,$ and $p_3$ are mutually orthogonal central projections, $p_1 M$ is a finite type $I$ von Neumann algebra, $p_2 M $ is a type $II_1$ von Neumann algebra, and $p_3 M$ is a properly infinite von Neumann algebra. Every element in $p_3 M$ writes as the sum of two commutators in $p_3 M$ (cf. \cite[Theorem 3.2 and Theorem 3.10]{FackHarpe1980sommes}). We can therefore conclude that $0=\eta p_3=p_3 \eta$.\smallskip

As in the proof of \Cref{thm C-star algebra being anti-derivable at zero}, by the Halving lemma, we can find two orthogonal equivalent projections $q_1$ and $q_2$ in $p_2 M$ such that $p_2 = q_1 + q_2$. Since $q_1-q_2$ has zero trace, and hence it can be written as the sum of ten commutators in $p_2 M$ \cite[Theorem 3.2 and Theorem 3.10]{FackHarpe1980sommes}, it follows that $\eta (q_1-q_2) =0$. We therefore have $0 =  \eta (q_1-q_2) (q_1-q_2) = \eta (q_1+q_2) = \eta p_2$.\smallskip

The von Neumann algebra $p_1 M$ decomposes as the direct sum of type $I_n$ von Neumann algebras with $n\in \mathbb{N}$. Suppose that one of these summands, whose unit is denoted by $z$, is a type $I_n$ von Neumann algebra with $n\geq 2$, and hence of the form $C(K, B(\ell_2^{n})),$ for some hyperstonean space $K.$ If $\xi$ and $\eta$ are two orthogonal unitary vectors in $\ell_2^{n}$, the constant function $u = \xi\otimes \eta$ has zero trace in $C(K, B(\ell_2^{n})),$ and thus the arguments above allow us to deduce that $\eta u = 0$, and then $0 = \eta u u^* = \eta (\xi\otimes \xi)$. The arbitrariness of $\xi$ and $\eta$ assures that $\eta z =0.$\smallskip

The above arguments show that $\eta$ belongs to the type $I_1$ part of $M$, as desired.\smallskip
  
For the reciprocal implication we assume that $(ii)$ is valid. Note that $M$ decomposes as an orthogonal sum of the form $M = p_{I_1} M \oplus^{\infty} N$, where $p_{I_1}$ is a central projection, $p_{I_1} M$ is the type $I_1$ part of $M$, and $\eta\in p_{I_1} M$. By observing that $p_{I_1} M$ is a commutative von Neumann algebra, given $a,b \in M$, we have $p_{I_1} [a,b] = [p_{I_1}a, p_{I_1}b ]=0.$ So, if we take $a,b \in M$ satisfying $a b =0$, we have
    \begin{align*}
        T(b) a +b T(a) &= b \eta a +b \eta a = 2 \eta ba = - 2 \eta [a,b] =-2 \eta p_{I_1} [a,b] =  0.
    \end{align*} So, $T$ is anti-derivable at zero.
\end{proof}

\textbf{Acknowledgements}\quad

\noindent First and third author supported by General Program of Shanghai Natural Science Foundation (Grant No.24ZR1415600).
Second author supported by MICIU/AEI/10.13039/501100011033 and ERDF/EU grant no. PID2021-122126NB-C31, by ``Maria de Maeztu'' Excellence Unit IMAG, reference CEX2020-001105-M, Junta de Andaluc{\'i}a grant FQM375, and by MOST (Ministry of Science and Technology of China) project no. G2023125007L. Third author also supported by ERDF/EU grant PID2021-122126NB-C31 and by China Scholarship Council grant No. 202306740016.


\subsection*{Data availability}

There is no data associate for the submission entitled ``New insights into linear maps {which are} anti-derivable at zero''.

\subsection*{Statements and Declarations}

The authors declare they have no financial nor conflicts of interests.


\bibliographystyle{plain}

\end{document}